
\magnification=1200
\voffset=0.3 true in
\parskip=2pt plus 1pt
\baselineskip=13.5pt plus 0.2pt minus 0.2pt
\lineskip=1.5pt plus0.5pt minus0.5pt
\lineskiplimit=0.5pt
\mathsurround=1pt
\font\tenmsx=msam10 \font\sevenmsx=msam7 \font\fivemsx=msam5
\newfam\msxfam
  \textfont\msxfam=\tenmsx \scriptfont\msxfam=\sevenmsx
\scriptscriptfont\msxfam=\fivemsx
\font\tenmsy=msbm10 \font\sevenmsy=msbm7 \font\fivemsy=msbm5
\newfam\msyfam
  \textfont\msyfam=\tenmsy \scriptfont\msyfam=\sevenmsy
\scriptscriptfont\msyfam=\fivemsy
\font\teneuf=eufm10 \font\seveneuf=eufm7 \font\fiveeuf=eufm5
\newfam\euffam
  \textfont\euffam=\teneuf \scriptfont\euffam=\seveneuf
\scriptscriptfont\euffam=\fiveeuf

\def\hexnumberA#1{\ifnum#1<10 \number#1\else
 \ifnum#1=10 A\else\ifnum#1=11 B\else\ifnum#1=12 C\else
 \ifnum#1=13 D\else\ifnum#1=14 E\else\ifnum#1=15 F\fi\fi\fi\fi\fi\fi\fi}
\def\MSX{\hexnumberA\msxfam}
\def\MSY{\hexnumberA\msyfam}

\def\text#1{\ifmmode
\mathchoice{\hbox{\tenrm #1}}
 {\hbox{\tenrm #1}}
 {\hbox{\sevenrm #1}}
 {\hbox{\fiverm #1}} \else\hbox{#1}\fi}
\def\matrix#1{\,\vcenter{\normalbaselines\mathsurround=0pt
  \ialign{\mathstrut\hfil$##$\hfil&&\quad\hfil$##$\hfil\crcr
  \mathstrut\crcr\noalign{\kern-\baselineskip}
  #1\crcr\mathstrut\crcr\noalign{\kern-\baselineskip}}}\,}
\def\fraction#1#2{{#1\over #2}}

\mathchardef\endproofmark="0\MSX03
\def\endofproof{\hphantom{\endproofmark}\hfill\llap{$\endproofmark$} }

\font\heading=cmr17 \font\footing=cmr9 \font\sectionfont=cmr10
\long\outer\def\Thm#1#2{\medbreak{\bf#1}{\sl#2}\medbreak}
\def\Prf#1{\medbreak{\bf#1}}

\def\Endrmk{\medbreak}
\def\Sct#1\par{\bigbreak\bigskip\bigskip\centerline{\sectionfont\S#1}
  \par\penalty10000\medskip}

\long\def\skipover#1\endskipover{}

\vglue1.5true in
\def\endproofmark{\hbox{\bf QED }}
\font\sectionfont=cmr12
\footline={\tenrm\hss --- \folio\ ---\hss}

\def\refindent{40pt}
\def\refhead{\sectionfont References}\def\reffont{\footing}
\def\Ref{\bigbreak\bigskip\bigskip\centerline{\refhead}\medskip
  \baselineskip=12pt plus.3pt minus.3pt\frenchspacing\reffont
  \everypar={\parindent=0pt\Refitem}\parindent=0pt
}

\def\Refitem [#1]#2:#3\par{\hangindent=\refindent\hangafter=1
  \rlap{[#1]}\hskip\hangindent #2: #3\par\smallbreak}
\bigskip
{\heading
\par\centerline {           Linear bound for abelian automorphisms}
\medskip
\par\centerline {                of varieties of general type}
}
\medskip
\par\centerline {                       (Version 1.0)}
\bigskip
{\footing
\par\centerline {                         Xiao, Gang}
\par\centerline {           Universit\'e de Nice, xiao@aurora.unice.fr}
}
\bigskip\bigskip\bigskip
    The aim of this paper is to prove the following
\medskip
    \Thm{Theorem 1.}{ Let $G$ be a finite abelian group acting faithfully on a
complex
smooth project variety $X$ of general type with numerically effective (nef)
canonical divisor, of dimension $n$. Then
                $$|G| \, \le  \, C(n)K_{X}^{n}  \
,$$
where $C(n)$ depends only on $n$.}
\bigskip
    We refer to the Introduction of [Ca-Sch] for a nice account of the history
for the study of bounds of automorphism groups of varieties of general type.
The authors of that paper have also shown a polynomial bound for abelian
automorphism groups.

    To prove Theorem 1, the only major obstacle to a generalisation of our
argument for surfaces [X] is the lack of a theorem of minimal models in higher
dimension: the basic idea is to find a pencil on $X$, whose general fibres are
invariant under the action of $G$, then use induction on $n$. To do so one
needs
bounded globally generatedness of pluricanonical sheaves, and vanishing
theorems. Unfortunately, these theorems currently exist only for varieties
with extra conditions which are not preserved by fibres. Therefore we consider
the problem for varieties in a more general category, as is done in [Ca-Sch].
Our main observation in Theorem 2 is that in the polynomial bound of Theorem
0.1 of [Ca-Sch], most copies of $d$ may be compensated by the ambient dimension
$N$, leading thus to a linear bound.

    All the coefficients in this paper are effective. As they are usually very
big and probably much bigger than the reality, we have preferred not to write
them down explicitely. Also no special effort has been made to optimise these
coefficients, so as to give a more readable presentation.
\bigskip
    \Sct 1. Preparation for induction\par
\medskip
    The main ingredient of the proof is the following theorem.
\medskip
    \Thm{Theorem 2.}{  Let $X\subset \mathchar"0\MSY50
_{\mathchar"0\MSY43 }^{N}$ be a non-degenerated irreducible variety of
dimension $n$ and degree $d$, with $\kappa (X)\ge
0$. Let $G$ be a finite abelian group acting
linearly on $\mathchar"0\MSY50 ^{N-1}$, leaving $X$ invariant. Then
                      $$|G| \, \le  \, c\left(n,d/N\right)d
\ ,$$
where the coefficient $c(n,d/N)$ depends only on $n$ and $d/N$, and is
increasing
with both variables.}
\medskip
    Theorem 1 follows directely from Theorem 2 and the following lemma, if we
send $X$ into $\mathchar"0\MSY50 ^{p_{m}(X)-1}$ by the linear system $|mK_
{X}|$ (take $d=m^{n}K_{X}^{n}$ and $N=p_{m}(X)$).
\medskip
    \Thm{Lemma 1.}{  For each positive integer $n$, there exist $M(n)\in
\mathchar"0\MSY4E $ and $a(n)\in \mathchar"0\MSY52
^{+}$,
such that:
\par     Let $X$ be a smooth projective variety of general type with nef
canonical
divisor, of dimension $n$. Then there is an integer $m\le
M(n)$, such that $|mK_{X}|$ is
free, $\Phi _{mK}$ is birational to its image, and
                          $$p_{m}(X)\ge a(n)m^{n}K_
{X}^{n}  \ .$$}
\par     \Prf{ Proof. }We follow the argument of [Wi], Theorem 1.1. By [Ko],
there is a
$r=2(n+2)(n+2)!$ such that $|mK_{X}|$ is free when $m\ge
r$. As $\Phi _{rK_{X}}$ is then generically
finite onto its image, we may take $n-1$ general divisors $D_
{1},\ldots ,D_{n-1}$ in $|rK_{X}|$
such that the intersection $D_{1}\cdots D_{n-1}$ is a smooth irreducible and
reduced curve
$C$.
\par     By successively applying the vanishing theorem of Kawamata-Viehweg, it
is
easy to see that the natural map
          $$H^{0}(X,((n-1)r+2)K_{X}) \, \longrightarrow
\, H^{0}(C,K_{C}+K_{X}|_{C})$$
is surjective. As $\deg(K_{C}+K_{X}|_{C})\ge 2g(C)+1$, $K_
{C}+K_{X}|_{C}$ is very ample on $C$, and
      $$p_{(n-1)r+4}(X)\ge h^{0}(C,K_{C}+K_{X}|_{C
})\ge {\fraction {1}{2}}\deg(K_{C}+K_{X}|_{C})={\fraction
{1}{2}}((n-1)r^{n}+2r^{n-1})K^{n}_{X}  \ .$$
This implies that $\Phi _{((n-1)r+2)K_{X}}$ is birational onto its image, so we
may take
            $$M(n)=(n-1)r+2  \ ,  \ a(n)={\fraction
{1}{2}}\left[{\fraction {   r}{(n-1)r+2}}\right]^{n-1
}  \ .$$
\endofproof\Endrmk
\bigskip
    \Sct 2. Proof of Theorem 2\par
\medskip
    First, we make the following remark which will be used in the argument.

    \Thm{Lemma 2.}{  In theorem 2, we may replace "irreducible variety" by
"subscheme whose components are of non-negative Kodaira dimension", under the
extra condition that each component of $X$ contains a point with trivial
stabiliser.}
\par     \Prf{ Proof. }Let $X_{1},\ldots ,X_{k}$ be the irreducible components
of $X$, $H_
{i}$ the minimal
subspace of $\mathchar"0\MSY50 ^{N-1}$ containing $X_
{i}$, with $d_{i}=\deg(X_{i})$, $N_{i}=\dim(H_{i})+1$. We have
$\sum _{i}d_{i}=d$, $\sum _{i}N_{i}\ge N$ as $X$ is non-degenerate, hence there
is an $i$ such that
$d_{i}/N_{i}\le d/N$. Let $G_{i}$ be the stabiliser of $X_
{i}$. Clearly $[G:G_{i}]\le d/d_{i}$, and $G_{i}$ acts
faithfully on $G_{i}$ by hypothesis. Thus we are reduced to the study of the
pair
$(X_{i}\subset H_{i},G_{i})$. \endofproof\Endrmk
\medskip
    Our basic tool is the induced linear actions of $G$ on the spaces $H^
{0}(\mathchar"024F _{X}(m))$
for $m\ge 1$. By the hypothesis that $G$ is finite abelian, such an action is
diagonalisable, hence $H^{0}(\mathchar"024F _{X}(m))$ has a basis cosisting of
semi-invariant
vectors. Each semi-invariant corresponds to a character of $G$.
\medskip
    \Prf{ Definition. }$H^{0}(\mathchar"024F _{X}(m))$ is called {\it uniquely
decomposable}, if
different semi-invariants correspond to different characters, or equivalently
if there are no more than $h^{0}(\mathchar"024F _{X
}(m))$ semi-invariants.\Endrmk
\medskip
    \Thm{Lemma 3.}{  We may assume that $H=H^{0}(\mathchar"024F
_{X}(1))$ is uniquely decomposable.}
\par     \Prf{ Proof. }We may obviously assume $N=\dim(H)$. Consider the
decomposition of
$H$ into eigenspaces $H=H_{1}\oplus \cdots \oplus
H_{k}$, with
                $$N_{1}\ge \ldots \ge N_{r}>1  \ ,
\ N_{r+1}=\cdots =N_{k}=1  \ ,$$
where $N_{i}=\dim(H_{i})$. For each $i\le r$, let $\pi
_{i}\colon  \, \mathchar"0\MSY50 ^{N-1}=\mathchar"0\MSY50
(H^{\vee }) \, \mathrel{\smash-\smash-\mathchar"0221
} \, \mathchar"0\MSY50 (H_{i}^{\vee })$ be the
projection with centre $\mathchar"0\MSY50 (H_{i}^{\vee
\perp})$, then let $T_{i}$ be a general hyperplane in
$\mathchar"0\MSY50 (H_{i}^{\vee })$, and $Y_{i}$ the moving part of the divisor
in $X$ cut out by $\pi
^{-1}_{i}(T_{i})$. Now fix
$i$ such that $Y_{i}$ is of maximal degree among these divisors.
\par     For an index $j\le r$ with $j\ne i$, let $P_
{i}$ be the minimal subspace of $\mathchar"0\MSY50
(H_{j}^{\vee })$
containing $\pi _{j}(Y_{i})$. Then $P_{j}$ is at least a hyperplane of
$\mathchar"0\MSY50
(H_{j}^{\vee })$, for otherwise
$\deg(Y_{j})>\deg(Y_{i})$ as one sees by taking a general hyperplane containing
$P_
{i}$,
taking into account that $Y_{i}$ being moving, no component of it can be
contained
in the fixed part of $\pi _{j}$.
\par     This means that there is at most one section in $H_
{j}$ vanishing on $Y_{i}$. Also,
no section in $H_{j}$ vanishes on $Y_{i}$ if $j>r$, for such a section is
proportionally
rigid, therefore it would vanish on a Zariski dense open subset of $X$,
contradicting the non-degeneratedness of $X$.
\par     Let $P$ be the minimal subspace of $\mathchar"0\MSY50
^{N-1}$ containing $Y_{i}$, with $N'=\dim(P)+1$. As
$Y_{i}$ is invariant under the action of $G$, $P$ is the intersection of $N-N'$
invariant
hyperplanes. So from the above, we get $N'\ge N-r\ge
N/2$. We also have
          $$n'=\dim(Y_{i})=n-1 \, , \, d'=\deg(Y_{i
})\le d \, ,$$
and $Y_{i}\subset P$ verifies the conditions of Lemma 2 from the easy addition
formula of
Kodaira dimensions. By induction,
              $$|G| \, \le  \, c(n',d'/N')d' \, \le
\, c(n-1,2d/N)d \, .$$
\endofproof\Endrmk
\bigskip
    We note by $D(X)$ the vector space ${ \text{div}
}(X)\otimes \mathchar"0\MSY51 $, where ${ \text{div}
}(X)$ is the additive
group of Weil divisors on $X$ (without taking linear equivalence). We are
interested in the finite subset $\Sigma \subset D(X)$ formed by the images of
the $N$ invariant
hyperplane sections, $H^{0}(\mathchar"024F _{X}(1))$ being uniquely
decomposable. For any positive
integer $m$ and $p_{1},\ldots ,p_{m}\in \Sigma $, the sum $p_
{1}+\cdots +p_{m}$ corresponds to an invariant
divisor in $|\mathchar"024F _{X}(m)|$, hence to a semi-invariant of $H^
{0}(\mathchar"024F _{X}(m))$. Our aim is to find
an $m$ such that semi-invariants of this kind outnumber $h^
{0}(\mathchar"024F _{X}(m))$, so that
$H^{0}(\mathchar"024F _{X}(m))$ is not uniquely decomposable.
\medskip
    \Prf{ Definition. }Let $\Sigma $ be a finite set in a vector space. The
{\it
dimension\/} of $\Sigma $ is the dimension of the convex hull of $\Sigma
$. Also, we define
                    $$m\Sigma  \, = \, \left\lbrace
\mathop{ \sum }\limits ^{ m}_{i=1}p_{i} \, | \, p_
{i}\in \Sigma \right\rbrace $$
to be the set of all sums of $m$ points in $\Sigma
$. The cardinal of $\Sigma $ will be denoted
by $|\Sigma |$.\Endrmk
\medskip
    \Thm{Lemma 4.}{  If $\Sigma $ is a finite set of dimension $\delta
$, then
              $$|m\Sigma | \, \ge  \, \left(\matrix
{m+\delta -1\cr \delta \cr }\right)|\Sigma |-\left(\matrix
{m+\delta -1\cr \delta -1\cr }\right)(m-1)$$
for any $m\ge 2$.}
\par     \Prf{ Proof. }Choose a point $p_{0}\in \Sigma
$ which is a summet of the convex hull of $\Sigma
$,
and let $\Sigma '=\Sigma \backslash \left\lbrace p_
{0}\right\rbrace $. When $|\Sigma |>\delta +1$, we may assume that $\Sigma
'$ is also of dimension
$\delta $. We can find $\delta $ points $p_{1},\ldots
,p_{\delta }\in \Sigma '$ such that the convex hull of
$\left\lbrace p_{0},p_{1},\ldots ,p_{\delta }\right\rbrace
$ is a simplex polyhedron whose intersection with the convex hull
of $\Sigma '$ is just the face generated by $p_{1},\ldots
,p_{\delta }$. In particular, all points of
the form $\mathop{ \sum }\limits ^{ \delta }_{i=0}n_
{i}p_{i}$ with $n_{0}>0$, $n_{i}\ge 0$, and $\mathop
{ \sum }\limits ^{ \delta }_{i=0}n_{i}=m$ are in $m\Sigma
$ but not $m\Sigma '$. As the
number of such points equals the number of $\delta
+1$-partitions of $m-1$, we get
                  $$|m\Sigma | \, \ge  \, |m\Sigma
'|+\left(\matrix{m+\delta -1\cr \delta \cr }\right)
\ .$$
The lemma follows by induction.  \endofproof\Endrmk
\medskip
    \Prf{ Remark. }The inequality of the lemma remains true when $\Sigma
$ is a set of
dimension $>\delta $: one has only to take a generic projection of $\Sigma
$ to a space of
dimension $\delta $.\Endrmk
\bigskip
    \Thm{Lemma 5.}{  Let $\Lambda $ be a linear system of affine dimension $N$
on a projective
variety $X$. Let $\left\lbrace D_{1},\ldots ,D_{k}\right\rbrace
$ be a set of generators of $\Lambda $, and $p_{i}$ the point in
$D(X)$ corresponding to $D_{i}$, for each $i$. Let $\delta
$ be the dimension of the set
$\left\lbrace p_{1},\ldots ,p_{k}\right\rbrace $. Then map $\Phi
_{\Lambda }\colon  \, X \, \mathrel{\smash-\smash-\mathchar"0221
} \, \mathchar"0\MSY50 ^{N-1}$ factors through a rational map
$\psi \colon  \, X \, \mathrel{\smash-\smash-\mathchar"0221
} \, \mathchar"0\MSY50 ^{\delta }$.}
\par     \Prf{ Proof. }We may obviously assume $k=N$. Let $A$ be the additive
subgroup of
$D(X)$ of rank $\delta $ generated by all the vectors $v_
{ij}=p_{i}-p_{j}$, and $\left\lbrace u_{1},\ldots
,u_{\delta }\right\rbrace $ a basis
of $A$. Let $f_{i}$ be the rational function on $X$ corresponding to $u_
{i}$. We have a
rational map $\psi \colon  \, X \, \mathrel{\smash-\smash-\mathchar"0221
} \, \mathchar"0\MSY50 ^{\delta }$ defined by $(1:f_
{1}:\cdots :f_{\delta })$.
\par     Now for any $i$ and $j$, we have $v_{ij}=\sum
_{k}n_{ijk}u_{k}$ with $n_{ijk}\in \mathchar"0\MSY5A
$, so if $g_{ij}$ is the
rational function corresponding to $v_{ij}$, we have $g_
{ij}=\mathop{ \prod }\limits ^{ \delta }_{k=1}f^{n_
{ijk}}_{k}$. $\Phi _{\Lambda }$ is defined,
say, by $\left(1:g_{21}:\cdots :g_{\delta 1}\right)$. Therefore if we define a
rational map
                          $$\alpha \colon  \, \mathchar"0\MSY50
^{\delta } \, \mathrel{\smash-\smash-\mathchar"0221
} \, \mathchar"0\MSY50 ^{N-1}$$
by
              $$\alpha (1:x_{1}:\cdots :x_{\delta
})=(1\colon \mathop{ \prod }\limits ^{ \delta }_{k=1
}x^{n_{11k}}_{k}\colon \mathop{ \prod }\limits ^{
\delta }_{k=1}x^{n_{21k}}_{k}:\cdots \colon \mathop
{ \prod }\limits ^{ \delta }_{k=1}x^{n_{\delta 1k}
}_{k})  \ ,$$
we have $\Phi _{\Lambda }=\alpha \circ \psi $ birationally. \endofproof\Endrmk
\medskip
    \Thm{Lemma 6.}{  Let $h_{m}(d,n)=\left(\matrix
{m+n-1\cr n\cr }\right)d+\left(\matrix{m+n-2\cr n-2\cr
}\right)$. For each positive integer $m$, we
have
                  $$h^{0}(\mathchar"024F _{X}(m))
\, \le  \, h_{m}(d,n) \text{.}$$}
\par     \Prf{ Proof. }This follows by double induction on $m$ and $n$, and the
exact
sequence
        $$0 \, \longrightarrow  \, H^{0}(\mathchar"024F
_{X}(m-1))  \ \longrightarrow  \, H^{0}(\mathchar"024F
_{X}(m)) \, \longrightarrow  \, H^{0}(\mathchar"024F
_{Y}(m))  \ ,$$
where $Y$ is a general hypereplane section of $X$. By additivity of Kodaira
dimensions, the induction on $n$ starts from curves of genus $\ge
1$. \endofproof\Endrmk
\medskip
    \Thm{Lemma 7.}{  Let $m$ be an integer greater than or equal to $(n+1)
{\fraction {  d}{N-n-1}}-n$. Then
$H^{0}(\mathchar"024F _{X}(m))$ is not uniquely decomposable under the induced
action of $G$.}
\par     \Prf{ Proof. }Let $\Sigma \in D(X)$ be the finite set corresponding to
invariant
hyperplane sections. As different points in $m\Sigma
$ correspond to non-proportional
semi-invariants in $H^{0}(\mathchar"024F _{X}(m))$, we have only to show
$|m\Sigma
|>h_{m}(d,n)$. But we have
$|m\Sigma |\ge \left(\matrix{m+n\cr n+1\cr }\right)N-\left(\matrix
{m+n\cr n\cr }\right)(m-1)$ by Lemma 4, as the dimension of $\Sigma
$ is at least $n+1$ due
to Lemma 5. The rest is straightforward computation. \endofproof\Endrmk
\medskip
    Now to use induction on $n$ to prove Theorem 2, choose an $m<(n+1)(n+2)
{\fraction {d}{N}}$ as
in Lemma 7, $F$ a general membre of the moving part of a pencil in $|mD|$
corresponding to a plane in $H^{0}(\mathchar"024F
_{X}(m))$ consisting of semi-invariants. $F$ is
invariant under the action of $G$. As $F$ is moving while $X$ has only finitely
many
invariant hyperplane sections due to Lemma 3, $F$ is non-degenerated in
$\mathchar"0\MSY50
^{N-1}$.
Therefore the couple $(F\subset \mathchar"0\MSY50
^{N-1},G)$ satisfies the conditions of Lemma 2, with
$\dim(F)=n-1$ and $\deg(F)\le \left[(n+1)(n+2){\fraction
{d}{N}}\right]d$, so
            $$|G| \, \le  \, \left[c\left(n-1,(n+1)(n+2)\left(
{\fraction {d}{N}}\right)^{2}\right)(n+1)(n+2){\fraction
{d}{N}}\right]d$$
by induction.
\medskip
    Finally, taking into account the proof of Lemma 3, Theorem 2 is shown by
setting
       $$c(n,{\fraction {d}{N}}) \, = \, \max\left\lbrace
c(n-1,2{\fraction {d}{N}}),c\left(n-1,(n+1)(n+2)\left(
{\fraction {d}{N}}\right)^{2}\right)(n+1)(n+2){\fraction
{d}{N}}\right\rbrace   \ .$$
\medskip
\Ref

    [Ca-Sch] Catanese, F. / Schneider, M.: Polynomial bounds for
abelian groups of automorphisms, to appear

    [Ko] Koll\'ar, J.: Effectif base point freeness, Math. Ann. 296, 595-605
(1993)

    [Wi] Wilson, P.: On complex algebraic varieties of general type, Symposia
Math. 24, 65-73 (1981)

    [X] Xiao, G.: On abelian automorphism group of a surface of general
type, Invent. Math. 102, 619-631 (1990)

\end